\documentclass[twocolumn,10pt]{asme2ej}

\usepackage{epsfig}
\usepackage{amssymb}
\usepackage{lineno}
\usepackage[tbtags]{amsmath}
\usepackage{graphicx,floatrow,threeparttable}
\title{Sliding control for single-degree-of-freedom fractional oscillators}

\author{Jian Yuan
    \affiliation{
	Institute of System Science \\and Mathematics\\
	Naval Aeronautical and \\Astronautical University\\
	Yantai, 264001, P.R.China\\
    Email: yuanjianscar@gmail.com
    }	
}

\author{Bao Shi
    \affiliation{
	Institute of System Science \\and Mathematics\\
	Naval Aeronautical and \\Astronautical University\\
	Yantai, 264001, P.R.China\\
    Email: baoshi781@sohu.com
    }	
}

\author{Guozhong Xiu\\
        \affiliation{
	Institute of System Science and Mathematics\\
	Naval Aeronautical and Astronautical University\\
	Yantai, 264001, P.R.China\\
    Email: xiuguozhong2013@163.com
    }	
}

\begin{document}

\maketitle

\begin{abstract}
{\it This paper proposes fractional sliding control designs for single-degree-of-freedom fractional oscillators respectively of the Kelvin-Voigt type, the modified Kelvin-Voigt type and D\"{u}ffing type, whose dynamical behaviors are described by second-order differential equations involving fractional derivatives. Firstly, the differential equations of motion are transformed into non-commensurate fractional state equations by introducing state variables with physical significance. Secondly, fractional sliding manifolds are constructed and stability of the corresponding sliding dynamics is addressed via the infinite state approach and Lyapunov stability theory. Thirdly, sliding control laws and adaptive sliding laws are designed for fractional oscillators respectively in cases that the bound of the external exciting force is known or unknown. Finally, numerical simulations are carried out to validate the above control designs.
}
\end{abstract}

\section{Introduction}
Viscoelastic materials are widely used to control or reduce vibrations and sound radiation in aerospace structures, industrial machines, civil engineering structures, etc [1]. In the quest for analyzing and synthesizing viscoelastically damped structures, one of the primary tasks is to model the rheological behavior of viscoelastic materials. However, this is also a challenge work due to the high dependency on viscoelastic materials of operational and environmental factors such as vibration frequency, external temperature, pre-loads, etc [2].

Recently, application of fractional calculus representations for describing the stress-strain relations for viscoelastic materials (shortly termed as fractional constitute relations) becomes a popular target. Over the conventional integer order constitutive relations, the fractional ones have vast superiority. The first attractive feature is that they are capable of fitting experimental results perfectly and describing mechanical properties accurately in both the frequency and time domain with only three to five empirical parameters [3]. The second is that they are not only consistent with the physical principles involved [4] and the molecule theory [5], but also represent the fading memory effect [3] and high energy dissipation capacity [6]. Finally, from mathematical perspectives the fractional constitutive equations and the resulting fractional differential equations of vibratory motion are compact and analytic [7]. Nowadays, via a large number of experiments, many types of fractional order constitutive relations have been developed, such as the fractional Kelvin-Voigt model [3], the fractional Zener model [4], and the fractional Pritz model [8], etc.

In engineering practice, mechanical models for various viscoelastically damped structures are simplified to SDOF or MDOF oscillators. On the other hand, investigations on typical oscillators provide theoretical basis for dealing with intricate mechanical systems. As a result, research on fractional oscillators becomes a promise subject. There are varieties of papers developing numerical approaches and methods for studying dynamics responses of fractional oscillators. Investigations on dynamical responses of SDOF linear and nonlinear fractional oscillators, MDOF fractional oscillators and infinite-DOF fractional oscillators have been reviewed in [9]. However, few results on theoretical analysis for fractional oscillators have been found in literature, which is partly attributed to imperfection in the theory of fractional differential equations, such as the initialization theory, state space representation and the stability theory. Asymptotically steady state behavior of fractional oscillators have been studied in [10, 11]. Based on functional analytic approach, criterias for existence and behavior of solutions have been obtained in [12-14], and particularly in which the impulsive response function for the linear SDOF fractional oscillator is derived. The asymptotically steady state response of fractional oscillators with more than one fractional derivatives have been analyzed in [15]. Taking the memory effect and prehistory into consideration, the historical effect and initialization problems for fractionally damped vibration equations have been proposed in [16-20].

Using fractional-order controllers is a suitable way to the more efficient control of fractional systems [21]. Various fractional control techniques have been proposed during the last two decades, such as CRONE control [22], fractional PID control [23], fractional sliding control [24, 25], fractional adaptive control [26, 27], fractional optimal control [28, 29], etc. However, control designs for fractional oscillators are very limited. In [30, 31], the second-order fractional differential equations of motion have been transformed into commensurate-ordered state equations in order to propose control designs for SDOF fractional oscillators. However, due to the fact that the orders of fractional derivatives in constitute relations for viscoelastic materials are usually real numbers lying in between zero and one, dimension of the resulting commensurate fractional state equations is normally too large to be analyzed theoretically. To overcome the above difficulty, in this study we introduce another types of state variables with physical significance to transform the fractional differential equations into non-commensurate fractional state equations. Based on this idea, dimensions of the corresponding fractional state equations are rather low, namely three for fractional Kelvin-Voigt oscillators and fractional D\"{u}ffing oscillators, four for fractional modified Kelvin-Voigt oscillators. Sliding control designs for the above three types of fractional oscillators are proposed based on infinite state approach and Lyapunov stability theory.

The rest of this paper is organized as follows. Section 2 gives some preliminaries on fractional calculus and differential equations of motion for SDOF fractional oscillators of Kelvin-Voigt type, modified Kelvin-Voigt type and D\"{u}ffing type. Section 3 proposes sliding control designs for the above three types of fractional oscillators. Section 4 carries out numerical simulations to verify the viability and effectiveness of the proposed control designs. Finally, the paper is concluded in section 5.

\section{Preliminaries}
\textbf{Definition 1} The Riemann-Liouville fractional integral for the function $f\left( t \right)$ is defined as
\begin{equation}
{}_{a}I_{t}^{\alpha }f\left( t \right)=\frac{1}{\Gamma \left( \alpha  \right)}\int_{a}^{t}{{{\left( t-\tau  \right)}^{\alpha -1}}f\left( \tau  \right)d\tau }
\end{equation}	
where $\alpha \in {{R}^{+}}$ is an non-integer order of the factional integral, the subscripts $a$ and $t$ are lower and upper terminals respectively.\\
\textbf{Definition 2} The Caputo definition of fractional derivatives is
\begin{equation}
{}_{a}D_{t}^{\alpha }f\left( t \right)=\frac{1}{\Gamma \left( n-\alpha  \right)}\int_{a}^{t}{\frac{{{f}^{\left( n \right)}}\left( \tau  \right)d\tau }{{{\left( t-\tau  \right)}^{\alpha -n+1}}}},n-1<\alpha <n
\end{equation}	
\textbf{Lemma 1} The frequency distributed model for the fractional integrator [32-34]
The input of the Riemann-Liouville integral is denoted by $v\left( t \right)$ and output $x\left( t \right)$, then ${}_{a}I_{t}^{\alpha }v(t)$ is equivalent to
\begin{equation}
\begin{cases}
\frac{\partial z(\omega ,t)}{\partial t}=-\omega z(\omega ,t)+v(t) \\
x(t)={}_{a}I_{t}^{\alpha }v\left( t \right)=\int_{0}^{+\infty }{{{\mu }_{\alpha }}(\omega )z(\omega ,t)d\omega } \\
\end{cases}
\end{equation}		
with ${{\mu }_{\alpha }}\left( \omega  \right)=\frac{\sin \left( \alpha \pi  \right)}{\pi }{{\omega }^{-\alpha }}$.\\
System (3) is the frequency distributed model for fractional integrator, which is also named as the diffusive representation.\\
The differential equation of motion for SDOF fractional Kelvin-Voigt oscillators is
\begin{equation}	
m\ddot{x}\left( t \right)+c{{D}^{\alpha }}x\left( t \right)+kx\left( t \right)=f\left( t \right)
\end{equation}	
where $m$ is the mass of the oscillator, $c$ is the coefficient of the fractional damping, $k$ is the stiffness coefficient, $\alpha \in \left( 0,1 \right)$ is the order of fractional derivative, $f\left( t \right)$ is external exciting force applied to the fractional oscillator.\\
The differential equation of motion for SDOF fractional modified Kelvin-Voigt oscillators is
\begin{equation}
m\ddot{x}\left( t \right)+{{c}_{1}}{{D}^{{{\alpha }_{1}}}}x\left( t \right)+{{c}_{2}}{{D}^{{{\alpha }_{2}}}}x\left( t \right)+kx\left( t \right)=f\left( t \right)
\end{equation}	
where $0<{{\alpha }_{1}}<{{\alpha }_{2}}<1$.\\
The differential equation of motion for fractional D\"{u}ffing oscillators is
\begin{equation}
\ddot{x}\left( t \right)+c{{D}^{\alpha }}x\left( t \right)\text{+}ax+b{{x}^{3}}\left( t \right)=f\left( t \right)
\end{equation}

\section{Sliding control design for SDOF fractional oscillators}
To propose sliding control design for fractional oscillators, the foremost task is to transform the fractional differential equations of motion (4) (5) or (6) into state equations. In contrast to integer-order derivatives, composition between fractional ones ${{D}^{\alpha }}{{D}^{\beta }}x\left( t \right)={{D}^{\alpha }}^{+\beta }x\left( t \right)={{D}^{\beta }}{{D}^{\alpha }}x\left( t \right)$ holds if and only if in zero initialization case: $x\left( t \right)\equiv 0$ for $t\le 0$. In physics, this means that fractional oscillators are at rest in prehistory. However, in nonzero initialization case, the above composition is more complicated. In this paper, zero initialization for fractional oscillators are assumed in transforming fractional differential equations into state equations. Nonzero initialization case will be investigated in our future work.

\subsection{Sliding control for SDOF fractional Kelvin-Voigt oscillators}
Letting ${{x}_{1}}=x$, ${{x}_{2}}={{D}^{\alpha }}x$, ${{x}_{3}}=\dot{x}$, the differential equation of motion (4) for the single degree of freedom fractional Kelvin-Voigt oscillators is transformed into the following state equations
\begin{equation}
\begin{cases}
{{D}^{\alpha }}{{x}_{1}}={{x}_{2}} \\
{{D}^{1-\alpha }}{{x}_{2}}={{x}_{3}} \\
{{{\dot{x}}}_{3}}=-\frac{k}{m}{{x}_{1}}-\frac{c}{m}{{x}_{2}}+\frac{f}{m}+\frac{u}{m} \\
\end{cases}
\end{equation}
In the following, we propose sliding control design for SDOF fractional oscillators of Kelvin-Voigt type (4) based on the above state space equations (7). Firstly, a fractional sliding surface is constructed to present the desired dynamics and secondly a switching control law is determined to verify sliding condition.\\
A fractional sliding surface is constructed as
\begin{equation}
s\left( t \right)={{x}_{3}}+\frac{k}{m}\int_{0}^{t}{{{x}_{1}}}\left( \tau  \right)d\tau +\frac{c}{m}\int_{0}^{t}{{{x}_{2}}}\left( \tau  \right)d\tau
\end{equation}	
Differentiating Eq.(8), leading to the following

$$\dot{s}={{x}_{3}}+\frac{k}{m}{{x}_{1}}+\frac{c}{m}{{x}_{2}}.$$
Letting $\dot{s}=0$, we derive the following sliding mode dynamics
\begin{equation}
\begin{cases}
{{D}^{\alpha }}{{x}_{1}}={{x}_{2}} \\
{{D}^{1-\alpha }}{{x}_{2}}={{x}_{3}} \\
{{{\dot{x}}}_{3}}=-\frac{k}{m}{{x}_{1}}-\frac{c}{m}{{x}_{2}} \\
\end{cases}
\end{equation}
Obviously, Eq. (9) are non-commensurate fractional differential equations. To analyze the stability of sliding dynamics, the infinite state approach based on the frequency distributed model is utilized.\\
Defining the following Lyapunov candidate
\[{{V}_{1}}\left( t \right)=\frac{1}{2}k{{x}_{1}}^{2}+\frac{c}{2}\int_{0}^{\infty }{{{\mu }_{1-\alpha }}\left( \omega  \right){{z}^{2}}\left( \omega ,t \right)d\omega }+\frac{1}{2}m{{\dot{x}}_{3}}^{2}\]
where $z\left( \omega ,t \right)$ is an infinite state variable in the following frequency distributed model for the Caputo derivative:
\begin{equation}
\begin{cases}
\frac{\partial z\left( \omega ,t \right)}{\partial t}=-\omega z\left( \omega ,t \right)+\dot{x} \\
{{D}^{\alpha }}x\left( t \right)=\int_{0}^{\infty }{{{\mu }_{1-\alpha }}\left( \omega  \right)z\left( \omega ,t \right)d\omega } \\
\end{cases}
\end{equation}
By taking the first order derivative of ${{V}_{1}}\left( t \right)$with respect to time yields
\[{{\dot{V}}_{1}}=-c\int_{0}^{\infty }{\omega {{\mu }_{1-\alpha }}}\left( \omega  \right){{z}^{2}}(\omega ,t)d\omega <0\]
which implies the stability of the sliding mode dynamics (9).\\
For the sake of the sliding control law, a second Lyapunov candidate is selected as
$${{V}_{2}}=\frac{1}{2}{{s}^{2}}$$.\\
By taking its derivative with respect to time yields
\begin{equation}
\begin{split}
{{\dot{V}}_{2}}&=s\dot{s}\\
&=s\left[ {{{\dot{x}}}_{3}}+\frac{k}{m}{{x}_{1}}+\frac{c}{m}{{x}_{2}} \right]
\end{split}
\end{equation}
Substituting the state equation (7) of the SDOF fractional oscillators into (11), it follows that
\begin{equation}
{{\dot{V}}_{2}}=s\left[ \frac{f\left( t \right)}{m}+\frac{u\left( t \right)}{m} \right]
\end{equation}
Firstly, the boundary of the external force $f\left( t \right)$is assumed to be known: $\left| f\left( t \right) \right|\le F$. In this case, Eq.(12) is calculated as
\begin{equation}
{{\dot{V}}_{2}}\le \frac{1}{m}\left[ \left| s \right|F+su \right]
\end{equation}
Based on (13), we construct the following control law
\begin{equation}
u\left( t \right)=-\left( F+\frac{{{\rho }_{1}}}{\sqrt{2}} \right)sign\left( s \right)-{{\rho }_{2}}s
\end{equation}
where ${{\rho }_{1}}$ and ${{\rho }_{2}}$ are positive constants, $sign\left( \cdot  \right)$ is the sign function.\\
Substituting (14) into (13), yields
\begin{equation*}
\begin{split}
{{\dot{V}}_{2}}&\le \frac{1}{m}\left[ \left| s \right|F-\left| s \right|\left( F+\frac{{{\rho }_{1}}}{\sqrt{2}} \right)-{{\rho }_{2}}{{s}^{2}} \right]\\
&=-\frac{{{\rho }_{1}}}{m}\sqrt{{{V}_{2}}}-2{{\rho }_{2}}{{V}_{2}}
\end{split}
\end{equation*}
In terms of Lyapunov stability theorem, we derive that $s\to 0$ and ${{x}_{1}}\left( t \right)\to 0$ as $t\to \infty $.\\
In the following, we propose the control design for the SDOF fractional oscillators in the case that the boundary of the external force is unknown. The adaptive control technique will be utilized to estimate the unknown boundary.  The estimation of $F$ is denoted as $\hat{F}\left( t \right)$.\\
Another Lyapunov candidate is chosen as
$${{\tilde{V}}_{2}}=\frac{1}{2}\left[ {{s}^{2}}+\frac{1}{\mu }{{\left( \hat{F}-F \right)}^{2}} \right]$$
where $\mu $ is a positive constant.\\
By taking its derivative with respect to time, one has
\begin{equation}
\begin{split}
{{\dot{\tilde{V}}}_{2}}&=s\dot{s}+\frac{1}{\mu }\left( \hat{F}-F \right)\dot{\hat{F}}\\
&=\frac{s}{m}\left[ f\left( t \right)+u\left( t \right) \right]+\frac{1}{\mu }\left( \hat{F}-F \right)\dot{\hat{F}}
\end{split}
\end{equation}
Based on (15), the following control law and adaptive law are constructed as
	\[u\left( t \right)=-\left( \hat{F}\left( t \right)+\frac{{{\rho }_{1}}}{\sqrt{2}} \right)sign\left( s \right)-{{\rho }_{2}}s\]
	\[\hat{F}\left( t \right)=\frac{\mu }{m}\left| s \right|\]	
Substituting (16) (17) into (15), and in terms of the relation $sign\left( s \right)s=\left| s \right|$, we derive that

\begin{equation}
\begin{split}
{{\dot{\tilde{V}}}_{2}}&=\frac{1}{m}\left[ f\left( t \right)s-\left( \hat{F}+{{\rho }_{1}} \right)\left| s \right|-{{\rho }_{2}}{{s}^{2}} \right]+\frac{1}{m}\left( \hat{F}-F \right)\left| s \right|\\
&=\frac{1}{m}\left[ f\left( t \right)s-\left( F+{{\rho }_{1}} \right)\left| s \right|-{{\rho }_{2}}{{s}^{2}} \right]
\end{split}
\end{equation}
Substituting the relation $f\left( t \right)s\le \left| f\left( t \right)s \right|\le \left| f\left( t \right) \right|\left| s \right|\le F\left| s \right|$ into (18), one has
$${{\dot{\tilde{V}}}_{2}}\le \frac{1}{m}\left[ -{{\rho }_{1}}\left| s \right|-{{\rho }_{2}}{{s}^{2}} \right]\le 0$$
In terms of Lyapunov stability theorem, we derive that $s\to 0$ and ${{x}_{1}}\left( t \right)\to 0$ as $t\to \infty $.
\subsection{Sliding control for SDOF fractional modified Kelvin-Voigt oscillators}
Letting ${{x}_{1}}=x$, ${{x}_{2}}={{D}^{{{\alpha }_{1}}}}x$, ${{x}_{3}}={{D}^{{{\alpha }_{2}}}}x$, ${{x}_{4}}=\dot{x}$, the differential equation of motion (5) for fractional oscillators with Kelvin-Voigt type is transformed into the following state equations
\begin{equation}
\begin{cases}
{{D}^{{{\alpha }_{1}}}}{{x}_{1}}={{x}_{2}} \\
{{D}^{{{\alpha }_{2}}-{{\alpha }_{1}}}}{{x}_{2}}={{x}_{3}} \\
{{D}^{{{\alpha }_{2}}}}{{x}_{3}}={{x}_{4}} \\
{{{\dot{x}}}_{4}}=-\frac{k}{m}{{x}_{1}}-\frac{{{c}_{1}}}{m}{{x}_{2}}-\frac{{{c}_{2}}}{m}{{x}_{3}}+\frac{f}{m}+\frac{u}{m}
\end{cases}
\end{equation}

A fractional sliding surface is constructed as
\begin{equation*}
\begin{split}
s\left( t \right)&={{x}_{4}}+\frac{k}{m}\int_{0}^{t}{{{x}_{1}}}\left( \tau  \right)d\tau +\frac{{{c}_{1}}}{m}\int_{0}^{t}{{{x}_{2}}}\left( \tau  \right)d\tau\\
&\quad +\frac{{{c}_{2}}}{m}\int_{0}^{t}{{{x}_{3}}}\left( \tau  \right)d\tau
\end{split}
\end{equation*}
Differentiating $s\left( t \right)$, leading to the following
$$\dot{s}\left( t \right)={{x}_{4}}+\frac{k}{m}{{x}_{1}}+\frac{{{c}_{1}}}{m}{{x}_{2}}+\frac{{{c}_{2}}}{m}{{x}_{3}}$$
Letting $\dot{s}=0$, yielding the following sliding mode dynamics
\begin{equation}
\begin{cases}
{{D}^{{{\alpha }_{1}}}}{{x}_{1}}={{x}_{2}} \\
{{D}^{{{\alpha }_{2}}-{{\alpha }_{1}}}}{{x}_{2}}={{x}_{3}} \\
{{D}^{1-{{\alpha }_{2}}}}{{x}_{3}}={{x}_{4}} \\
{{{\dot{x}}}_{4}}=-\frac{k}{m}{{x}_{1}}-\frac{c}{m}{{x}_{2}}-\frac{{{c}_{2}}}{m}{{x}_{3}} \\
\end{cases}
\end{equation}
To analyze the stability of sliding dynamics (20), the following Lyapunov candidate is defined

\begin{equation}
\begin{split}
{{V}_{3}}&=\frac{1}{2}kx_{1}^{2}+\frac{{{c}_{1}}}{2}\int_{0}^{\infty }{{{\mu }_{1-{{\alpha }_{1}}}}}\left( \omega  \right){{z}_{1}}(\omega ,t)d\omega\\
& \quad +\frac{{{c}_{2}}}{2}\int_{0}^{\infty }{{{\mu }_{1-{{\alpha }_{2}}}}}\left( \omega  \right){{z}_{2}}(\omega ,t)d\omega +\frac{1}{2}mx_{4}^{2}
\end{split}
\end{equation}
where ${{z}_{1}}\left( \omega ,t \right)$ and ${{z}_{2}}\left( \omega ,t \right)$ are infinite state variables in the following frequency distributed model for the Caputo derivatives
\[{{D}^{{{\alpha }_{1}}}}x={{I}^{1-{{\alpha }_{1}}}}\dot{x}\Leftrightarrow \left\{ \begin{aligned}
  & \frac{\partial {{z}_{1}}\left( \omega ,t \right)}{\partial t}=-\omega {{z}_{1}}\left( \omega ,t \right)+\dot{x} \\
 & {{D}^{{{\alpha }_{1}}}}x=\int_{0}^{\infty }{{{\mu }_{1-{{\alpha }_{1}}}}\left( \omega  \right){{z}_{1}}\left( \omega ,t \right)d\omega } \\
\end{aligned} \right.\]
\[{{D}^{{{\alpha }_{2}}}}x={{I}^{1-{{\alpha }_{2}}}}\dot{x}\Leftrightarrow \left\{ \begin{aligned}
  & \frac{\partial {{z}_{2}}\left( \omega ,t \right)}{\partial t}=-\omega {{z}_{1}}\left( \omega ,t \right)+\dot{x} \\
 & {{D}^{{{\alpha }_{2}}}}x=\int_{0}^{\infty }{{{\mu }_{1-{{\alpha }_{2}}}}\left( \omega  \right){{z}_{2}}\left( \omega ,t \right)d\omega } \\
\end{aligned} \right.\]
By taking the first order derivative of ${{V}_{3}}\left( t \right)$with respect to time yields
\begin{equation*}
\begin{split}
{{\dot{V}}_{3}}&=-{{c}_{1}}\int_{0}^{\infty }{\omega {{\mu }_{1-{{\alpha }_{1}}}}}\left( \omega  \right){{z}_{1}}^{2}(\omega ,t)d\omega \\
&\quad-{{c}_{2}}\int_{0}^{\infty }{\omega {{\mu }_{1-{{\alpha }_{2}}}}}\left( \omega  \right){{z}_{2}}^{2}(\omega ,t)d\omega \\
 &\quad<0
\end{split}
\end{equation*}
which implies the stability of the sliding mode dynamics (20).\\
For the sake of the sliding control law, the following Lyapunov candidate is selected as
\[{{V}_{4}}=\frac{1}{2}{{s}^{2}}\]
By taking its derivative with respect to time yields

\begin{equation}
{{\dot{V}}_{4}}=s\dot{s}=s\left[ {{x}_{4}}+\frac{k}{m}{{x}_{1}}+\frac{{{c}_{1}}}{m}{{x}_{2}}+\frac{{{c}_{2}}}{m}{{x}_{3}} \right]
\end{equation}
Substituting the state equation (19) of the SDOF fractional oscillators into (21)£¬one has

\begin{equation}
{{\dot{V}}_{4}}=s\left[ \frac{f}{m}+\frac{u}{m} \right]
\end{equation}
In the case that the boundary of the external force is known: $\left| f\left( t \right) \right|\le F$, Eq. (22) is then calculated as

\begin{equation}
{{\dot{V}}_{4}}\le \frac{1}{m}\left[ \left| s \right|F+su \right]
\end{equation}
Based on (23), we construct the following control law

\begin{equation}
u\left( t \right)=-\left( F+\frac{{{\rho }_{1}}}{\sqrt{2}} \right)sign\left( s \right)-{{\rho }_{2}}s
\end{equation}
Substituting (24) into (23), one has
\begin{equation}
\begin{split}
{{{\dot{V}}}_{4}}&\le \frac{1}{m}\left[ \left| s \right|F-\left| s \right|\left( F+\frac{{{\rho }_{1}}}{\sqrt{2}} \right)-{{\rho }_{2}}{{s}^{2}} \right] \\
 & \quad=-\frac{{{\rho }_{1}}}{m}\sqrt{{{V}_{4}}}-2{{\rho }_{2}}{{V}_{4}}
\end{split}
\end{equation}
In terms of Lyapunov stability theorem, we derive that $s\to 0$ and ${{x}_{1}}\left( t \right)\to 0$ as $t\to \infty $.\\
In the case that the boundary of the external force is unknown, the same control law and adaptive law as (16) and (17) can be designed utilizing the adaptive technique.
\subsection{Sliding control for fractional D\"{u}ffing oscillators}
Letting ${{x}_{1}}=x$, ${{x}_{2}}={{D}^{\alpha }}x$, ${{x}_{3}}=\dot{x}$, the differential equation of motion (6) for the Fractional Duffing oscillators is transformed into the following state equations
\begin{equation}
\left\{ \begin{aligned}
  & {{D}^{\alpha }}{{x}_{1}}={{x}_{2}} \\
 & {{D}^{1-\alpha }}{{x}_{2}}={{x}_{3}} \\
 & {{{\dot{x}}}_{3}}=\text{-}a{{x}_{1}}-bx_{1}^{3}-c{{x}_{2}}+f+u \\
\end{aligned} \right.
\end{equation}
A fractional sliding surface is constructed as
\[s={{x}_{3}}+k\int_{0}^{t}{{{x}_{1}}}\left( \tau  \right)d\tau +c\int_{0}^{t}{{{x}_{2}}}\left( \tau  \right)d\tau \]
Differentiating $s\left( t \right)$, leading to the following
\[\dot{s}={{\dot{x}}_{3}}+k{{x}_{1}}+c{{x}_{2}}\]
Letting $\dot{s}=0$, we obtain the following sliding mode dynamics
\begin{equation}
\left\{ \begin{aligned}
  & {{D}^{\alpha }}{{x}_{1}}={{x}_{2}} \\
 & {{D}^{1-\alpha }}{{x}_{2}}={{x}_{3}} \\
 & {{{\dot{x}}}_{3}}=-k{{x}_{1}}-c{{x}_{2}} \\
\end{aligned} \right.
\end{equation}
To analyze the stability of sliding dynamics (26), the following Lyapunov candidate is defined
\[{{V}_{5}}=\frac{1}{2}kx_{1}^{2}+\frac{1}{2}x_{3}^{2}+\frac{c}{2}\int_{0}^{\infty }{{{\mu }_{1-\alpha }}}\left( \omega  \right){{z}^{2}}(\omega ,t)d\omega \]
where $z\left( \omega ,t \right)$ is an infinite state variable in the frequency distributed model (10).\\
By taking its derivative with respect to time yields
\begin{equation}
\begin{split}
{{\dot{V}}_{5}}&=k{{x}_{1}}{{x}_{3}}+{{x}_{3}}{{\dot{x}}_{3}}+c{{x}_{3}}{{x}_{2}}-c\int_{0}^{\infty }{\omega {{\mu }_{1-\alpha }}}\left( \omega  \right){{z}^{2}}d\omega\\
&={{x}_{3}}\left[ {{{\dot{x}}}_{3}}+c{{x}_{2}}+k{{x}_{1}} \right]-c\int_{0}^{\infty }{\omega {{\mu }_{1-\alpha }}}\left( \omega  \right){{z}^{2}}d\omega
\end{split}
\end{equation}
Substituting the state equation (26) into (27), one has
\[{{\dot{V}}_{5}}=-c\int_{0}^{\infty }{\omega {{\mu }_{1-\alpha }}}\left( \omega  \right){{z}^{2}}d\omega <0\]
which implies the stability of the sliding mode dynamics (26).\\
For the sake of sliding control law, the following Lyapunov candidate is selected as
\[{{V}_{6}}\left( t \right)=\frac{1}{2}{{s}^{2}}\]
By taking its derivative with respect to time yields
\begin{equation}
{{\dot{V}}_{6}}\left( t \right)=s\dot{s}=s\left[ {{{\dot{x}}}_{3}}+k{{x}_{1}}+c{{x}_{2}} \right]
\end{equation}
Substituting the state equation (26) into (28), one has
\begin{equation}
{{\dot{V}}_{6}}\left( t \right)=s\left[ \left( \text{-}a+k \right){{x}_{1}}-bx_{1}^{3}+f+u \right]
\end{equation}
We construct the following control law
\begin{equation}
u\left( t \right)=bx_{1}^{3}\text{+}\left( a-k \right){{x}_{1}}+v\left( t \right)
\end{equation}
Substituting (30) into (29) yields
\begin{equation}
{{\dot{V}}_{6}}\left( t \right)=s\left[ f+v \right]=sf+sv\le \left| s \right|F+sv
\end{equation}

We select $v\left( t \right)=-\left( F+\frac{{{\rho }_{1}}}{\sqrt{2}} \right)sign\left( s \right)-{{\rho }_{s}}s$ and substitute it into (31), then Eq.(31) follows that
\begin{equation*}
\begin{split}
{{{\dot{V}}}_{6}}&\le \frac{1}{m}\left[ \left| s \right|F-\left| s \right|\left( F+\frac{{{\rho }_{1}}}{\sqrt{2}} \right)-{{\rho }_{2}}{{s}^{2}} \right] \\
 &=-\frac{{{\rho }_{1}}}{m}\sqrt{{{V}_{6}}}-2{{\rho }_{2}}{{V}_{6}} \\
\end{split}
\end{equation*}
In terms of Lyapunov stability theorem, we derive that $s\to 0$ and ${{x}_{1}}\left( t \right)\to 0$ as $t\to \infty $.\\
In the case that the boundary of the external force is unknown, the same control law and adaptive law as (16) and (17) can be designed utilizing the adaptive technique.

\section{Numerical simulations}
In this section, we present numerical simulations in MATLAB to evaluate the performance of the sliding mode control for SDOF fractional oscillators.

Parameters in the fractional oscillator (4) are taken respectively as $m=1$, $c=0.4$, $k=2$, $\alpha =0.56$, the external force is assumed to be $f\left( t \right)=30\cos 6t$. The forced vibration response of the fractional oscillator (4) is shown in Fig.1.
\begin{figure}
\centering
\includegraphics[scale=0.6]{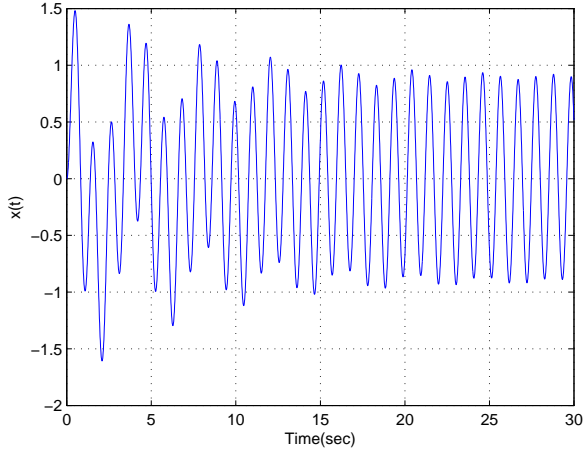}
\caption{Vibration response of the SDOF fractional Kelvin-Voigt oscillators}
\end{figure}

In the sliding control law (14), parameters are taken as $F+\frac{{{\rho }_{1}}}{\sqrt{2}}=35$, ${{\rho }_{2}}=2$. The performance is illustrated in Fig.2. In the sliding control law (16) and adaptive law (17), parameters are taken as $\frac{{{\rho }_{1}}}{\sqrt{2}}=25$, ${{\rho }_{2}}=5$, $\mu =5$. The performance is illustrated in Fig.3.

\begin{figure}
\centering
\includegraphics[scale=0.6]{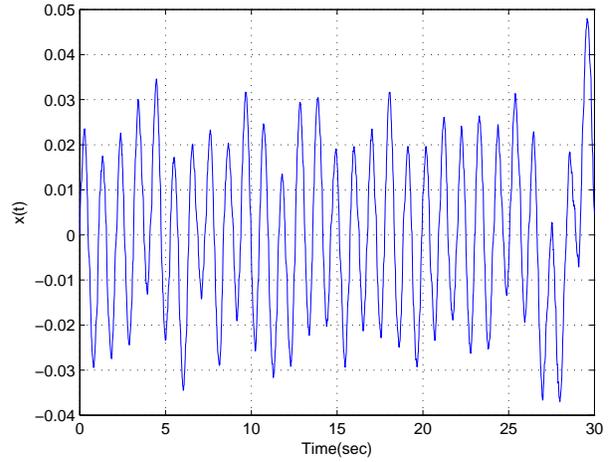}
\caption{State response of the sliding control for SDOF fractional Kelvin-Voigt oscillators}
\end{figure}
\begin{figure}
\centering
\includegraphics[scale=0.6]{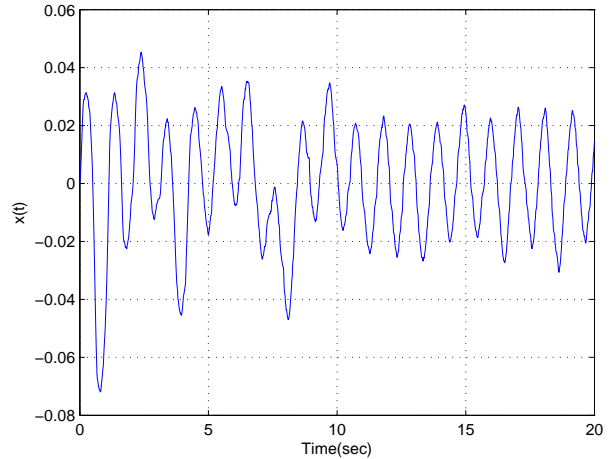}
\caption{State response of the adaptive sliding control for SDOF fractional Kelvin-Voigt oscillators}
\end{figure}
In the numerical simulations for the fractional oscillator with modified Kelvin-Voigt type, the mass are taken as  $m=1$, the stiffness coefficient $k=2$, the fractional damping coefficient ${{c}_{1}}=0.4$ and ${{c}_{2}}=0.2$, the orders of the fractional derivatives ${{\alpha }_{1}}=0.56$ and ${{\alpha }_{2}}=0.2$, the external force $f\left( t \right)=30\cos 6t$.
The forced vibration response of the fractional oscillator (5) is shown in Fig.4. In the sliding control law (24), parameters are taken as $F+\frac{{{\rho }_{1}}}{\sqrt{2}}=31$, ${{\rho }_{2}}=1$. The performance is illustrated in Fig.5.
\begin{figure}
\centering
\includegraphics[scale=0.6]{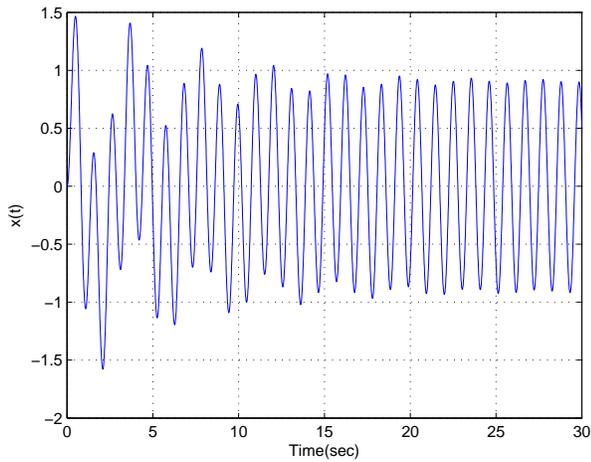}
\caption{Vibration response of the SDOF fractional modified Kelvin-Voigt oscillators}
\end{figure}

\begin{figure}
\centering
\includegraphics[scale=0.6]{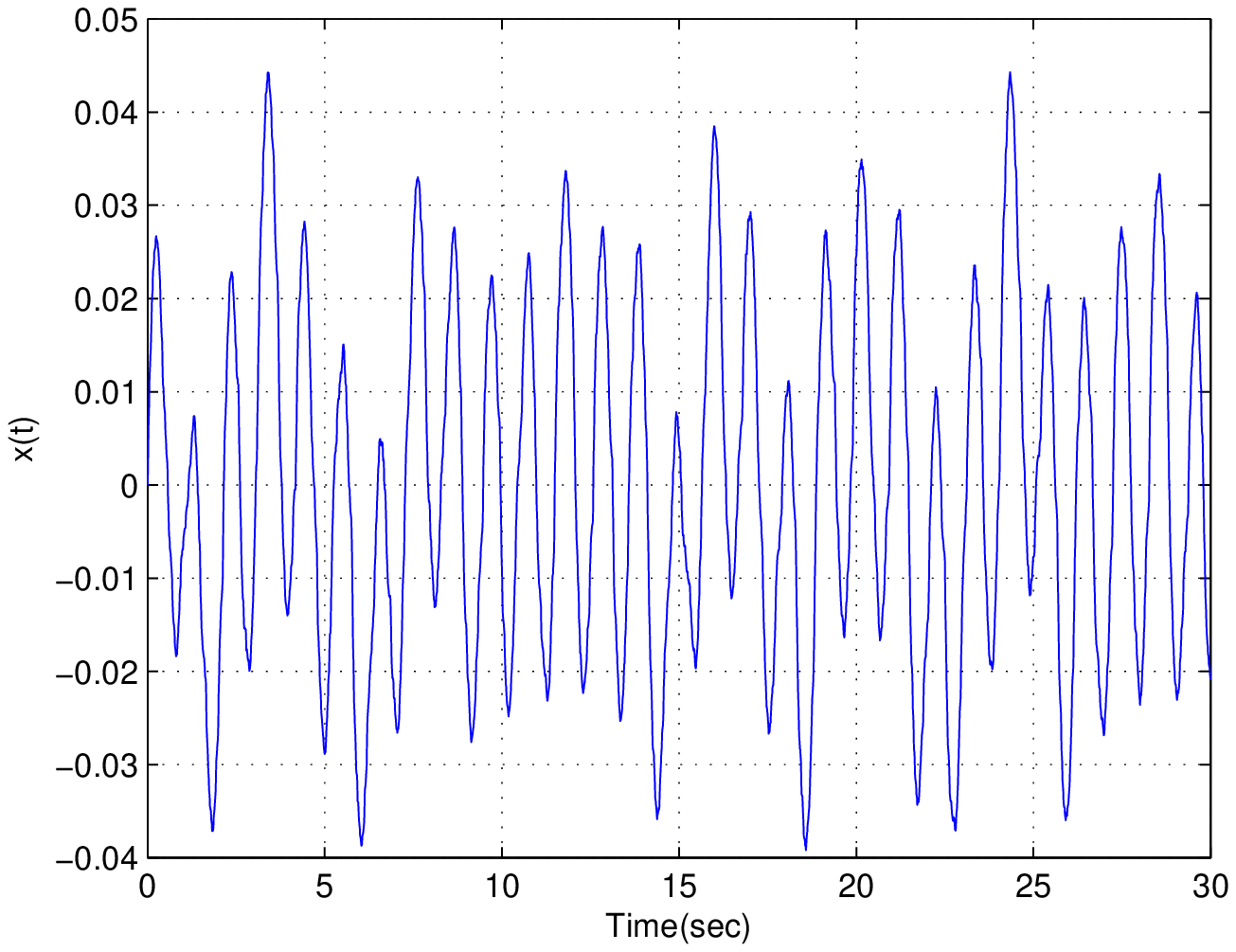}
\caption{State response of the sliding control for SDOF fractional modified Kelvin-Voigt oscillators}
\end{figure}
\begin{figure}
\centering
\includegraphics[scale=0.6]{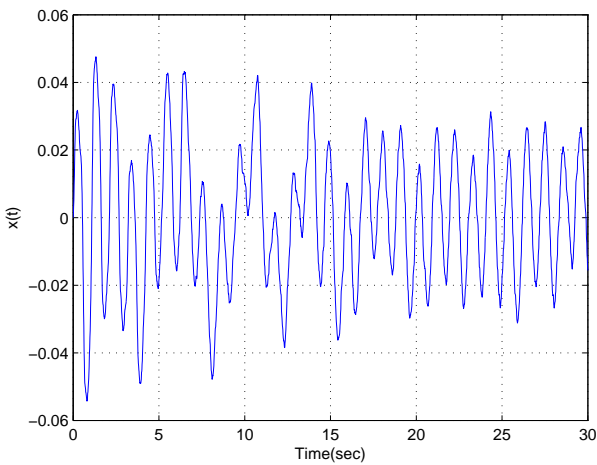}
\caption{State response of the adaptive sliding control for SDOF fractional modified Kelvin-Voigt oscillators}
\end{figure}

In the numerical simulations for the fractional D\"{u}ffing oscillators (6), parameters are taken as $a=-2$, $b=4$, $c=0.4$, $\alpha =0.56$, the external force $f\left( t \right)=30\cos 6t$. The forced vibration response of the fractional D\"{u}ffing oscillators (6) is shown in Fig.7.

In the sliding control law (30), parameters are taken as $k=1$, $F+\frac{{{\rho }_{1}}}{\sqrt{2}}=31$£¬${{\rho }_{2}}=2$. The performance is illustrated in Fig.8.
\begin{figure}
\centering
\includegraphics[scale=0.6]{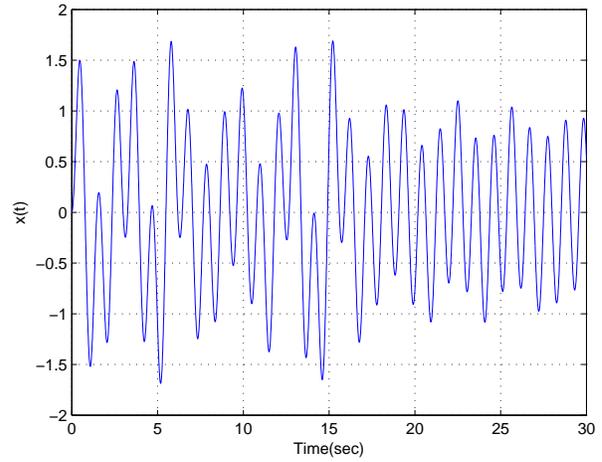}
\caption{Vibration response of the fractional D\"{u}ffing oscillators}
\end{figure}

\begin{figure}
\centering
\includegraphics[scale=0.6]{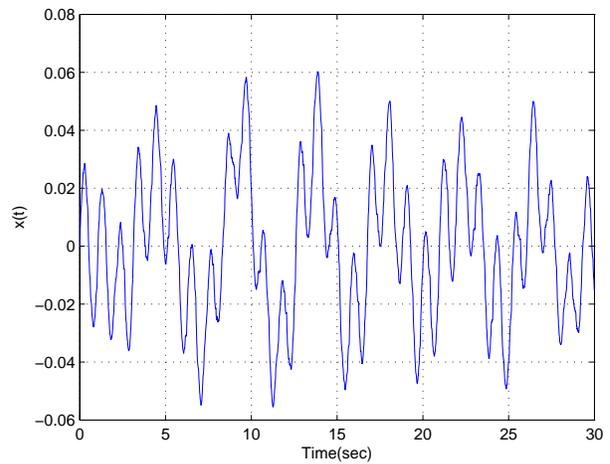}
\caption{State response of the sliding control for the fractional D\"{u}ffing oscillators}
\end{figure}

\section{Conclusions}
This paper have proposed fractional sliding control designs for SDOF fractional oscillators of the Kelvin-Voigt type, the modified Kelvin-Voigt type and D\"{u}ffing type. Differential equations of motion have been transformed into non-commensurate fractional state equations by introducing state variables with physical significance. Fractional sliding manifolds have been constructed and stability of the corresponding sliding dynamics has been addressed via the infinite state approach and Lyapunov stability theory. Sliding control laws and adaptive sliding laws have been designed respectively for fractional oscillators in cases that the bound of the external exciting force is known or unknown. Finally, the viability and effectiveness of the above control designs have been validated by numerical simulations.

\begin{acknowledgment}
The author Yuan Jian expresses his thanks to Prof. Dong Kehai from Naval Aeronautical and Astronautical University, and Prof. Jiang Jianping from national University of Defense technology. All the authors acknowledge the valuable suggestions from the peer reviewers. This work was supported by the Natural Science Foundation of the Province Shandong of China titled Controls for fractional systems with applications to hypersonic vehicles (Grant Nos. ZR2014AM006).
\end{acknowledgment}
\textbf{References}\\
1.	S. Park, Analytical modeling of viscoelastic dampers for structural and vibration control, International Journal of Solids and Structures, 38 (2001) 8065-8092.\\
2.	A.M.G.d. Lima, N. Bouhaddi, D.A. Rade, M. Belonsi, A time-domain finite element model reduction method for viscoelastic linear and nonlinear systems, Latin American Journal of Solids and Structures, 12 (2015) 1182-1201.\\
3.	R.L. Bagley, Applications of Generalized Derivatives to Viscoelasticity, in, Air Force Institute of technology, 1979, pp. 133.\\
4.	R.L. Bagley, J. TORVIK, Fractional calculus-a different approach to the analysis of viscoelastically damped structures, AIAA journal, 21 (1983) 741-748.\\
5.	R.L. Bagley, P. Torvik, A theoretical basis for the application of fractional calculus to viscoelasticity, Journal of Rheology (1978-present), 27 (1983) 201-210.\\
6.	Z.-D. Xu, C. Xu, J. Hu, Equivalent fractional Kelvin model and experimental study on viscoelastic damper, J. Vib. Control, (2013) 1077546313513604.\\
7.	X. Moreau, C. Ramus-Serment, A. Oustaloup, Fractional differentiation in passive vibration control, Nonlinear Dynam., 29 (2002) 343-362.\\
8. T. Pritz, Five-parameter fractional derivative model for polymeric damping materials, Journal of Sound and Vibration, 265 (2003) 935-952.\\
9.	Y.A. Rossikhin, M.V. Shitikova, Application of fractional calculus for dynamic problems of solid mechanics: novel trends and recent results, Applied Mechanics Reviews, 63 (2010) 010801.\\
10.	J. Padovan, S. Chung, Y.H. Guo, Asymptotic steady state behavior of fractionally damped systems, J. Franklin. I., 324 (1987) 491-511.\\
11.	J. Padovan, Y. Guo, General response of viscoelastic systems modelled by fractional operators, J. Franklin. I., 325 (1988) 247-275.\\
12.	H. Beyer, S. Kempfle, Definition of physically consistent damping laws with fractional derivatives, ZAMM Journal of Applied Mathematics and Mechanics/Zeitschrift f\"{u}r Angewandte Mathematik und Mechanik, 75 (1995) 623-635.\\
13.	S. Kempfle, I. Sch?fer, H. Beyer, Fractional calculus via functional calculus: theory and applications, Nonlinear Dynam., 29 (2002) 99-127.\\
14.	I. Sch\"{a}fer, S. Kempfle, Impulse responses of fractional damped systems, Nonlinear Dynam., 38 (2004) 61-68.\\
15.	Y.A. Rossikhin, M.V. Shitikova, Analysis of rheological equations involving more than one fractional parameters by the use of the simplest mechanical systems based on these equations, Mechanics of Time-Dependent Materials, 5 (2001) 131-175.\\
16.	M. Fukunaga, On initial value problems of fractional differential equations, International Journal of Applied Mathematics, 9 (2002) 219-236.\\
17.	M. Fukunaga, On uniqueness of the solutions of initial value problems of ordinary fractional differential equations, International Journal of Applied Mathematics, 10 (2002) 177-190.\\
18. M. Fukunaga, A difference method for initial value problems for ordinary fractional differential equations II, International Journal of Applied Mathematics, 11 (2003) 215-244.\\
19. T.T. Hartley, C.F. Lorenzo, Control of initialized fractional-order systems, NASA Technical Report, (2002).\\
20.	C.F. Lorenzo, T.T. Hartley, Initialization of Fractional-Order Operators and Fractional Differential Equations, Journal of Computational and Nonlinear Dynamics, 3 (2008) 021101.\\
21.	I. Podlubny, Fractional differential equations: an introduction to fractional derivatives, fractional differential equations, to methods of their solution and some of their applications, Academic press, 1999.\\
22.	A. Oustaloup, X. Moreau, M. Nouillant, The CRONE suspension, Control Engineering Practice, 4 (1996) 1101-1108.\\
23.	I. Podlubny, Fractional-order systems and PI/sup/spl lambda//D/sup/spl mu//-controllers, Automatic Control, IEEE Transactions on, 44 (1999) 208-214.\\
24.	B.S. Jian Yuan, Wenqiang Ji, Adaptive sliding mode control of a novel class of fractional chaotic systems, Advances in Mathematical Physics, (2013).\\
25.	B. Bandyopadhyay, S. Kamal, Stabilization and Control of Fractional Order Systems: A Sliding Mode Approach, Springer, 2015.\\
26.	S. Ladaci, A. Charef, On Fractional Adaptive Control, Nonlinear Dynam., 43 (2006) 365-378.\\
27.	B. Shi, J. Yuan, C. Dong, On fractional Model Reference Adaptive Control, The Scientific World Journal, 2014 (2014) 521625.\\
28.	O.P. Agrawal, A general formulation and solution scheme for fractional optimal control problems, Nonlinear Dynam., 38 (2004) 323-337.\\
29.	D. Baleanu, J.A.T. Machado, A.C.J. Luo, A Formulation and Numerical Scheme for Fractional Optimal Control of Cylindrical Structures Subjected to General Initial Conditions, (2012).\\
30.	R.L. Bagley, R.A. Calico, Fractional order state equations for the control of viscoelastically damped structures, Journal of Guidance, Control, and Dynamics, 14 (1991) 304-311.\\
31.	A. Fenander, Modal synthesis when modeling damping by use of fractional derivatives, AIAA Journal, 34 (1996) 1051-1058.\\
32. G. Montseny, Diffusive representation of pseudo-differential time-operators, LAAS, 1998.\\
33.	J.C. Trigeassou, N. Maamri, J. Sabatier, A. Oustaloup, Transients of fractional-order integrator and derivatives, Signal, Image and Video Processing, 6 (2012) 359-372.\\
34.	J.C. Trigeassou, N. Maamri, J. Sabatier, A. Oustaloup, State variables and transients of fractional order differential systems, Comput. Math. Appl., 64 (2012) 3117-3140.\\

\end{document}